\documentclass[11pt]{amsart}
\usepackage{verbatim}
\usepackage{amsmath}
\usepackage{pdfsync}

\newtheorem{theorem}{Theorem}[section]

\newtheorem{lemma}{Lemma}[section]

\numberwithin{equation}{section}

\usepackage{hyperref}
\hypersetup{colorlinks,citecolor=blue,plainpages=false,hypertexnames=false}

\begin{document}

\title[Criteria for the ampleness of certain vector bundles]{Criteria for the ampleness of certain vector bundles}

\author{Indranil Biswas}

\address{Department of Mathematics, Shiv Nadar University, NH91, Tehsil
Dadri, Greater Noida, Uttar Pradesh 201314, India}

\email{indranil.biswas@snu.edu.in, indranil29@gmail.com}

\author{Vamsi Pritham Pingali}

\address{Department of Mathematics, Indian Institute of Science, Bangalore 560012, India}

\email{vamsipingali@iisc.ac.in}



\begin{abstract}
We prove that certain vector bundles over surfaces are ample if they are so when restricted to divisors, certain numerical criteria
hold, and they are semistable (with respect to $\det(E)$). This result is a higher-rank version of a theorem
of Schneider and Tancredi for vector bundles of rank two over surfaces. We also provide
counterexamples indicating that our theorem is sharp.
\end{abstract}

\maketitle

\section{Introduction}\label{sec:Intro}

A holomorphic vector bundle $E$ over a compact complex manifold $X$ is said to be ample if
$\mathcal{O}_{\mathbb{P}(E^{*})}(1)$ is ample over $\mathbb{P}(E^{*})$ (the projective bundle
over $X$ parametrising the lines in the fibres of $E^{*}$). Ample vector bundles
play an important role in algebraic geometry because of the various vanishing theorems that
ensue from ampleness. Many of these vanishing theorems have numerous geometric consequences.
It is therefore of interest to find criteria for ampleness.

For any line bundle $L$ over a projective manifold $X$, the Nakai-Moizeshon criterion
gives a numerical condition to decide the ampleness of $L$ \cite{NM}. It says that
$L$ is ample if for all $1\,\leq\, k\, \leq\, \dim X$, we have $c_1(L)^k.Y\,>\,0$ for every closed
subvariety $Y$ of dimension $k$ of $X$. However, no
numerical criterion can exist for deciding the ampleness of vector bundles \cite{FL}.
Notwithstanding this negative result, Schneider and Tancredi proved the following criterion
(that is not purely numerical) for rank two vector bundles over surfaces.

\begin{theorem}[{\cite[p.~134]{ScTa}}]
Let $E$ be a holomorphic vector bundle of rank two over a compact complex surface $X$. Assume
that $c_1(E)\,>\,0$ and that $E$ is semistable with respect to $\det(E)$. Suppose
$E\big\vert_C$ is ample for every closed curve $C\, \subset\, X$, and
$$(c_1(E)^2-2c_2(E)).X\,>\, 0,\ \ c_2(E).X\,>\,0.$$ Then $E$ is ample.
\label{thm:ST}
\end{theorem}

This result and its improvements carried out in \cite{ScTa2} are useful in studying the
ampleness of cotangent bundles \cite{Sc}. They also find use in an approach of Demailly to the
Green-Griffiths-Lang conjecture \cite{Dem}. In this paper we aim to generalise this theorem
to vector bundles of other ranks for further possible applications.

We were motivated by the following result of L\"ubke \cite{Lub}.

\begin{theorem}[{\cite[p.~313, Theorem 2.1]{Lub}}]
Let $(E,\, h)$ be a holomorphic Hermitian rank-$r$ vector bundle over a compact K\"ahler manifold $(X,\,
\omega)$ of dimension $n$. Suppose $F_h \wedge \omega^{n-1}\,=\,-\sqrt{-1}\lambda \omega^n$, where $F_h$
is the curvature of the Chern connection of $h$ and $\lambda\,>\,0$ is a constant. Assume
that $$c_1(E,\,h)\,=\,\frac{r\lambda}{2\pi} \omega .$$ Also, suppose there exists a positive
function $\psi$ such that either of the following holds:
\begin{enumerate}
\item $n\,=\,2$ and $c_1^2(E,h)-\frac{2r(r-1)}{r^2-2r+2}c_2(E,\,h)\,=\,\psi \omega^2$, or
\item $r\,=\,2$ and $c_1^2(E,h)-\frac{4(n-1)^2}{n^2-2n+2}c_2(E,\,h)\,=\,\psi \omega^2$.
\end{enumerate}
Then $h$ is Griffiths-positively curved, i.e., $\langle v,\,\sqrt{-1}F_h v\rangle$ is
a K\"ahler form whenever $v\,\neq\, 0$ is a vector in $E$.
\label{thm:Lub}
\end{theorem}

\indent The following is the main result of this paper.

\begin{theorem}\label{thm:Mainthm}
Let $E$ be a holomorphic vector bundle of rank $r$ over a compact complex manifold $X$ of
dimension two. Suppose $c_1(E)\,>\,0$ and $E$ is semistable with respect
to $\det(E)$. Also assume that $E$ restricted to every curve is ample, and that
that $(c_1^2-c_2)(E).X\,>\,0$. Then $E$ is ample if
\begin{equation}\label{gc}
\left(c_1^2(E)-\frac{2r(r-1)}{r^2-2r+2}c_2(E)\right).X\,>\,0.
\end{equation}
\end{theorem}

The proof of Theorem \ref{thm:Mainthm} (carried out in Section \ref{sec:proof}) uses the
existence of approximately Hermitian-Einstein metrics on semistable vector bundles \cite{Kob}.

In Section \ref{sec:counter} we provide examples to indicate that L\"ubke's Chern class inequality in
Theorem \ref{thm:Lub} cannot be dispensed with for $n\,=\,2$ (and arbitrary $r$).

\section{Proof of Theorem \ref{thm:Mainthm}}\label{sec:proof}

In this section, we prove Theorem \ref{thm:Mainthm}. The Nakai-Moishezon
criterion will be used \cite{NM}. Our aim is to show that $(c_1(\mathcal{O}_{\mathbb{P}(E^{*})}(1)))^d.Y\,>\,0$
for every subvariety $Y$ of $\mathbb{P}(E^{*})$ of dimension $d$. Let
$$
\pi\, \,:\,\, \mathbb{P}(E^{*})\, \longrightarrow\, X
$$
be the natural projection. If $\pi(Y)$ is a point, we are done trivially. If $\pi(Y)$ is a curve, then since $E$
restricted to $\pi(Y)$ is ample, we are done. So assume that $\pi(Y)\,=\, X$.

In this case we shall compute the intersection number by choosing
an appropriate smooth metric $h$ on $E$ and considering the Chern-Weil representative of the induced metric
$\widetilde{h}$ on $\mathcal{O}_{\mathbb{P}(E^{*})}(1)$. Firstly, we fix a K\"ahler form $\omega$ on $X$. Let
$$S\,\, \subset\,\, X$$
be the smallest Zariski closed proper subset such that $Y\setminus \pi^{-1}(S)$ consists of regular points of the
projection map over $X\setminus S$, i.e., $Y\setminus \pi^{-1}(S)$ is a smooth fibre bundle over $X\setminus S$.

For every $\epsilon$, there exists an approximate Hermitian-Einstein metric $h_{\epsilon}$ (with curvature
$F_{\epsilon}$) satisfying $c_1(h_{\epsilon})\,=\,\omega$ and \eqref{eq:approximateHE} \cite{Kob}. We shall
choose $\epsilon$ later. Let $\Theta_{\epsilon}\,=\,\frac{\sqrt{-1}}{2\pi}F_{\epsilon}$. Let $(p,\,[v])\,\in\,
Y$. The key point is that if we choose a holomorphic normal trivialisation of $E$ near $p$, then
\begin{gather}
c_1(\mathcal{O}_{\mathbb{P}(E^{*})}(1),\,\widetilde{h}_{\epsilon})(p,\,[v])\,=
\,\frac{\langle v,\, \pi^*\Theta_{\epsilon} v\rangle}{\langle v,\, v\rangle}+\omega_{FS},
\label{eq:normaltrivformula}
\end{gather}
where $\omega_{FS}$ is the Fubini-Study metric on the fibres of $\pi$ (restricted to $Y$). Therefore,
$$
c_1(\mathcal{O}_{\mathbb{P}(E^{*})}(1),\,\widetilde{h}_{\epsilon})^d(p,\,[v])
\,=\,
$$
$$
\frac{d(d-1)}{2}\left(\pi^*\frac{\langle v,\, \Theta_{\epsilon},\, v\rangle}{\langle v,
\,v\rangle}\right)^2\omega^{d-2}_{FS} + \omega_{FS}^d + d \pi^*\frac{\langle v,\, \Theta_{\epsilon},\, v\rangle}{\langle v,
\,v\rangle} \omega_{FS}^{d-1}
$$
\begin{equation}\label{eq:highestpowerinnormaltriv}
=\,\frac{d(d-1)}{2}\left(\pi^*\frac{\langle v,\, \Theta_{\epsilon},\, v\rangle}{\langle v,\,
v\rangle}\right)^2c_1(\mathcal{O}_{\mathbb{P}(E^{*})}(1),\,\widetilde{h}_{\epsilon})^{d-2}(p,\,[v]),
\end{equation}
where we noted that on a surface, at most two powers of $\frac{\langle v,\, \Theta_{\epsilon},\, v\rangle}{\langle v,\,
v\rangle}$ are non-zero, and since the dimension of the fibre is $d-2$,
we have $\omega_{FS}^{d-1}\,=\,0$ on $Y \cap \pi^{-1}(p)$ for any $p\,\in\, S$. Note that the expression
$$\frac{d(d-1)}{2}\left(\pi^*\frac{\langle v,\, \Theta_{\epsilon},\, v\rangle}{\langle v,\,
v\rangle}\right)^2c_1(\mathcal{O}_{\mathbb{P}(E^{*})}(1),\,\widetilde{h}_{\epsilon})^{d-2}(p,\,[v])$$
is independent of the choice of the local trivialisation of $E$. 

At this juncture, the following lemma, which is a pointwise assertion, will be used to lower bound the
right-hand-side of \eqref{eq:highestpowerinnormaltriv}.

\begin{lemma}\label{lem:pointwiseinequality}
For every $\epsilon \,>\,0$, let $\Theta_{\epsilon}$ be a (normalized) Chern curvature endomorphism of a Hermitian holomorphic
vector bundle $(E,\, h_{\epsilon})$, of rank $r$, at a point $p$ on a surface $X$. Let $v\,\in\, E_p$. Suppose $\omega$ is a K\"ahler form
at $p$ and $c_1(h_{\epsilon})\,=\,\omega$. Moreover, assume that for every $\epsilon\,>\,0$, there exists a trace-free endomorphism
$B_{\epsilon}$ of $E$ at $p$ satisfying $\big\vert (B_{\epsilon})_i^j \big\vert \,\leq\, \epsilon$ for all $i,\, j$, and
\begin{gather}
\Theta_{\epsilon}\wedge \omega \,=\, \frac{1}{r} \omega^2 + \frac{B_{\epsilon}}{2} \omega^2.
\label{eq:approximateHE}
\end{gather}
Then
\begin{gather}
c_1^2(h_{\epsilon})-\frac{2r(r-1)}{r^2-2r+2}c_2 (h_{\epsilon})\,\leq\, \frac{r^2}{r^2-2r+2}\left(\frac{\langle v,\, \Theta_{\epsilon} v
\rangle_{h_{\epsilon}}}{\langle
v,\, v\rangle_{h_{\epsilon}}}\right)^2+\frac{4r +r(r^2-1)\epsilon}{4(r^2-2r+2)}\epsilon\omega^2.
\label{ineq:pointwiseinequality}
\end{gather}
\end{lemma}

\begin{proof}
Since the statement does not depend on the choice of trivialisation as well as the choice of coordinates, we can assume that $h_{\epsilon}$
is given by the identity matrix at this point $p$, in other words, $h_{\epsilon}$ is the trivial
Hermitian structure. Moreover, without loss of generality, we assume that $v\,=\,(1,\,0,\,0)$. For ease of notation, we drop the $\epsilon$ subscript on $\Theta_{\epsilon}$. Lastly, we can choose coordinates so that
$$\omega\,\,=\,\,\displaystyle \sum_{i=1}^{r} \Theta_{i}^{i} \,\,=\,\, \sqrt{-1} dz^1 \wedge d\overline{z}^1+\sqrt{-1}dz^2\wedge d\overline{z}^2,$$ and
$\Theta_1^1 \,=\, \sqrt{-1} \mu_1 dz^1 \wedge d\overline{z}^1+\sqrt{-1}\mu_2 dz^2\wedge d\overline{z}^2$. By the approximate
Hermitian-Einstein condition \eqref{eq:approximateHE} we see that
\begin{gather}
\mu_1+\mu_2\,=\,\frac{2}{r}+(B_{\epsilon})_1^1\nonumber\\
\Theta_2^2+\Theta_3^3+\ldots\,=\,
(1-\mu_1)\sqrt{-1}dz^1\wedge d\overline{z}^1+\left(1-\frac{2}{r}+\mu_1-(B_{\epsilon})_1^1\right)
\sqrt{-1}dz^2\wedge d\overline{z}^2,\nonumber\\
\Theta_i^i\wedge \omega \,=\,\frac{1}{r} \omega^2 +\frac{1}{2}(B_{\epsilon})_i^i \omega^2,\nonumber\\
(\Theta_i^j)_{1\overline{1}}+(\Theta_i^j)_{2\overline{2}}\,=\,(B_{\epsilon})_i^j.\label{eq:deductionsfromapproxandcoordinates}
\end{gather}
Now,
\begin{gather}
\frac{c_1(h)^2-\frac{2r(r-1)}{r^2-2r+2}c_2(h)}{\sqrt{-1}dz^1 \wedge d\overline{z}^1 \wedge
\sqrt{-1}dz^2\wedge d\overline{z}^2} \,=\,
2 - \frac{r(r-1)}{r^2-2r+2}\frac{\sum_{i\neq j}\left( -\Theta_i ^j \overline{\Theta}_i^j+
\Theta_i^i \Theta_j^j \right)}{\sqrt{-1}dz^1 \wedge d\overline{z}^1 \wedge
\sqrt{-1}dz^2\wedge d\overline{z}^2} 
\nonumber \\
\leq\, 2 - \frac{r(r-1)}{r^2-2r+2}\times\nonumber \\
\frac{\sum_{i\neq j}\left( -((\Theta_i^j)_{1\overline{1}} \overline{(\Theta_i^j)_{2\overline{2}}}+
(\Theta_i^j)_{2\overline{2}} \overline{(\Theta_i^j)_{1\overline{1}}} )
\sqrt{-1}dz^1 \wedge d\overline{z}^1 \wedge \sqrt{-1}dz^2\wedge d\overline{z}^2+
\Theta_i^i \Theta_j^j \right)}{\sqrt{-1}dz^1 \wedge d\overline{z}^1 \wedge \sqrt{-1}dz^2\wedge d\overline{z}^2}.
\nonumber
\end{gather}
 Using \eqref{eq:deductionsfromapproxandcoordinates} we see that
\begin{gather}
\frac{c_1(h)^2-\frac{2r(r-1)}{r^2-2r+2}c_2(h)}{\sqrt{-1}dz^1 \wedge d\overline{z}^1 \wedge
\sqrt{-1}dz^2\wedge d\overline{z}^2}\,\leq \nonumber \\
2+\frac{r(r-1)\big\vert (B_{\epsilon})_i^j\big\vert^2}{r^2-2r+2} - \frac{r(r-1)\displaystyle \sum_{i\neq j}
\left(\Theta_i^i \Theta_j^j\right)}{(r^2-2r+2)\sqrt{-1}dz^1
\wedge d\overline{z}^1 \wedge \sqrt{-1}dz^2\wedge d\overline{z}^2} \nonumber \\
=2+\frac{r(r-1)\big\vert (B_{\epsilon})_i^j\big\vert^2}{r^2-2r+2} - \frac{r(r-1)\Bigg(\omega^2 -\sum_i (\Theta_i ^i)^2\Bigg)}{(r^2-2r+2)\sqrt{-1}dz^1
\wedge d\overline{z}^1 \wedge \sqrt{-1}dz^2\wedge d\overline{z}^2} \nonumber \\
\leq 2+\frac{r(r-1)\epsilon^2}{r^2-2r+2} -\frac{2r(r-1)}{r^2-2r+2} + \frac{r(r-1)\sum_i (\Theta_i ^i)^2}{(r^2-2r+2)\sqrt{-1}dz^1
\wedge d\overline{z}^1 \wedge \sqrt{-1}dz^2\wedge d\overline{z}^2} \nonumber \\
=\frac{4-2r+r(r-1)\epsilon^2}{r^2-2r+2}+\frac{r(r-1)\sum_i (\Theta_i ^i)^2}{(r^2-2r+2)\sqrt{-1}dz^1
\wedge d\overline{z}^1 \wedge \sqrt{-1}dz^2\wedge d\overline{z}^2}\nonumber \\
\leq \frac{4-2r+r(r-1)\epsilon^2}{r^2-2r+2}+\frac{r(r-1)(\Theta_1^1)^2}{(r^2-2r+2)\sqrt{-1}dz^1
\wedge d\overline{z}^1 \wedge \sqrt{-1}dz^2\wedge d\overline{z}^2}\nonumber \\
+\frac{2r(r-1)\sum_{i\geq 2} (\Theta_i ^i)_{1\bar{1}}(\Theta_i ^i)_{2\bar{2}}}{r^2-2r+2} \nonumber \\
=\frac{4-2r+r(r-1)\epsilon^2}{r^2-2r+2}+\frac{r(r-1)(\Theta_1^1)^2}{(r^2-2r+2)\sqrt{-1}dz^1
\wedge d\overline{z}^1 \wedge \sqrt{-1}dz^2\wedge d\overline{z}^2}\nonumber \\
+\frac{2r(r-1)\sum_{i\geq 2} (\Theta_i ^i)_{1\bar{1}}\left(\frac{2}{r}+(B_{\epsilon})^i_i- (\Theta_i ^i)_{1\bar{1}}\right)}{r^2-2r+2}\nonumber \\
\Rightarrow\,\, \frac{(r^2-2r+2)c_1(h)^2-2r(r-1)c_2(h)}{\sqrt{-1}dz^1 \wedge d\overline{z}^1 \wedge
\sqrt{-1}dz^2\wedge d\overline{z}^2}\,\, \leq\nonumber \\
4-2r+r(r-1)\epsilon^2+\frac{r(r-1)(\Theta_1^1)^2}{\sqrt{-1}dz^1
\wedge d\overline{z}^1 \wedge \sqrt{-1}dz^2\wedge d\overline{z}^2}\nonumber \\
+2r(r-1)\sum_{i\geq 2} (\Theta_i ^i)_{1\bar{1}}\left(\frac{2}{r}+(B_{\epsilon})^i_i- (\Theta_i ^i)_{1\bar{1}}\right).
\label{ineq:firstmajor}
\end{gather}
Now we want to maximize $f\,=\,\sum_{i\geq 2} (\Theta_i ^i)_{1\bar{1}}\left(\frac{2}{r}+(B_{\epsilon})^i_i-
(\Theta_i ^i)_{1\bar{1}}\right)$ subject to the condition $\sum_{i\geq 2} (\Theta_i^i)_{1\bar{1}}\,=\,1-\mu_1$. Clearly, $f$ tends to
$-\infty$ at infinity. Therefore, using Lagrange's multipliers we conclude that the maximum of $f$ is as follows
(we replace $\mu_1$ by $\mu$ for the remainder of this section):
\begin{gather}
\displaystyle \sum_{i\geq 2} \left(\left(\frac{1-\mu}{r-1}+\frac{(B_{\epsilon})_1^1}{2(r-1)}+\frac{(B_{\epsilon})_i^i}{2} \right) \left(\frac{2}{r}+(B_{\epsilon})_i^i \right) -\left(\frac{1-\mu}{r-1}+\frac{(B_{\epsilon})_1^1}{2(r-1)}+\frac{(B_{\epsilon})_i^i}{2} \right)^2 \right).
\end{gather}
Thus we have
\begin{gather}
\frac{(r^2-2r+2)c_1(h)^2-2r(r-1)c_2(h)}{\sqrt{-1}dz^1 \wedge d\overline{z}^1 \wedge
\sqrt{-1}dz^2\wedge d\overline{z}^2}\nonumber \\
\leq\,4-2r+r(r-1)\epsilon^2+\frac{r(r-1)(\Theta_1^1)^2}{\sqrt{-1}dz^1
\wedge d\overline{z}^1 \wedge \sqrt{-1}dz^2\wedge d\overline{z}^2}
\nonumber \\
+2r-4-2r(B_{\epsilon})_1^1-\frac{r}{2}((B_{\epsilon})_1^1)^2+\frac{r(r-1)\sum_{i\geq 2}((B_{\epsilon})_i^i)^2}{2} +2r\mu\left(\frac{2}{r}+(B_{\epsilon})_1^1-\mu \right) \nonumber \\
\leq\,\, \frac{r(r^2-1)\epsilon^2}{2}+2r\epsilon+\frac{r^2(\Theta_1^1)^2}{\sqrt{-1}dz^1
\wedge d\overline{z}^1 \wedge \sqrt{-1}dz^2\wedge d\overline{z}^2}.
\end{gather}
This completes the proof of the lemma.
\end{proof}

Returning to the situation at hand, choose $$\epsilon\,\,=\,\, \min \left(1,\,\, \frac{2\int ((r^2-2r+2)c_1(h)^2-2r(r-1)c_2(h))}{r(r^2+1)\int \omega^2} \right)
$$ and let $h\,=\,h_{\epsilon}$. Since $c_1(\mathcal{O}_{\mathbb{P}(E^{*})}(1),\,\widetilde{h}_{\epsilon})$
is a closed form representing the cohomology class $c_1(\mathcal{O}_{\mathbb{P}(E^{*})}(1))$, we have
$$
(c_1(\mathcal{O}_{\mathbb{P}(E^{*})}(1)))^d.Y\,\,=\,\,
\int_Y c_1(\mathcal{O}_{\mathbb{P}(E^{*})}(1),\widetilde{h}_{\epsilon})^d.
$$
Therefore, Lemma \ref{lem:pointwiseinequality} shows that
$$
(c_1(\mathcal{O}_{\mathbb{P}(E^{*})}(1)))^d.Y\,\,\geq\,\,
$$
$$
\frac{d(d-1)}{4r^2} \int_S ((r^2-2r+2)c_1^2-2r(r-1)c_2)(p) \left(\int_{\pi^{-1}(p)\cap Y} \omega_{FS}^{d-2}\right) .
$$
Now $\int_{\pi^{-1}(p)\cap Y} \omega_{FS}^{d-2}$ is the degree of $\pi^{-1}(p)\bigcap Y\,\subset\, \pi^{-1}(p)\,=\,
\mathbb{P}(E^{*}_p)$, and hence the given condition \eqref{gc} in Theorem \ref{thm:Mainthm} implies that
$$
(c_1(\mathcal{O}_{\mathbb{P}(E^{*})}(1)))^d.Y\,\,> \,\, 0.
$$
This completes the proof.

\section{Counterexamples}\label{sec:counter}

Consider vector bundles of rank $r$ on surfaces. In this section we provide counterexamples to show that if L\"ubke's condition
$$c_1^2(E).X\,\,>\,\,\frac{2r(r-1)}{r^2-2r+2}c_2(E).X$$ is not met, then the conclusion of the theorem cannot hold in general.
Likewise for semistability with respect to $\det(E)$. For $r\,=\,2$, the
counterexample for semistability was provided in \cite{ScTa} and the sharpness of the Chern class inequality was shown
in \cite{Su}. For the remainder of this section we assume that $r\,\geq\, 3$.

Indeed, just as in \cite{ScTa}, let $M$ be a Riemann surface of genus at least two, and let $F$ be a stable vector bundle of rank
two on $M$ such that $c_1(F)\,=\,0$ and the symmetric product $S^m F$ is stable for all $m\,\geq\, 1$.
Set $X\,=\,\mathbb{P}(F)$, and let $L\,=\,\mathcal{O}_{\mathbb{P}(F)}(1)$. (Then $c_1(L)^2\,=\,0$ and $L\big\vert_C$ is positive
for all curves $C$ on $X$.) Let $H$ be an ample line bundle on $X$, and
set $$E\,=\,L\oplus H \oplus H \oplus \ldots \oplus H .$$ Note that $E$ is not ample because its
quotient $L$ is not ample. However, $E$ is ample when restricted to curves. Then we have
\begin{gather}
c_1(E)\,=\,c_1(L)+(r-1)c_1(H),\nonumber \\
c_2(E)\,=\,(r-1)c_1(L)c_1(H)+\frac{(r-1)(r-2)}{2}c_1^2(H).\nonumber
\end{gather}
The slope of $L$ is $c_1(L).c_1(E)\,=\,(r-1)c_1(L)c_1(H)$ and that of $H$ is $$c_1(H).c_1(E)\,=\,c_1(H).c_1(L)+(r-1)c_1(H)^2.$$
Therefore, the semistability of $E$ (with respect to $\det(E)$) holds if and only if
\begin{gather}
(r-2)c_1(L).c_1(H)\,=\,(r-1)c_1(H)^2.
\nonumber
\end{gather}
Certainly this inequality cannot be met for an arbitrary $H$ (and thus ampleness can fail if semistability does not hold). Suppose
that semistability is met. Then we have the following:
\begin{gather}
c_1(E)^2.X\,=\,(r-1)^2c_1(H)^2.X+2(r-1)c_1(L)c_1(H).X=(r-1)rc_1(L)c_1(H).X, \nonumber \\
\frac{2r(r-1)}{r^2-2r+2}c_2(E).X\,=\,\frac{2r(r-1)}{r^2-2r+2}\frac{2(r-1)+(r-2)^2}{2} c_1(L)c_1(H)=c_1(E)^2.\nonumber
\end{gather}
Therefore, if L\"ubke's Chern class inequality is not met, ampleness cannot hold in general.

\section*{Acknowledgements}

The authors are sincerely grateful to the anonymous referee for detailed and useful feedback that enabled us to correct several crucial
errors. The second author is partially supported by the DST FIST program - 2021 [TPN - 700661], and a MATRICS grant
MTR/2020/000100 from SERB (Govt. of India). The first author is partially supported by a J. C. Bose Fellowship.

\end{document}